\documentclass[10pt]{amsart}
\pdfoutput=1
\usepackage{amsmath}
\usepackage{amsfonts}
\usepackage{amssymb}
\usepackage{amscd}
\usepackage{graphics}
\usepackage{graphicx}
\usepackage {xypic}
 \input{xy}
  \xyoption{all}

\newtheorem{theorem}{Theorem}[section]
\newtheorem{definition}[theorem]{Definition}
\newtheorem{proposition}[theorem]{Proposition}
\newtheorem{lemma}[theorem]{Lemma}
\newtheorem{corollary}[theorem]{Corollary}

\newtheorem{conjecture}[theorem]{Conjecture}

\title{The relative Burnside kernel - the elementary abelian case}
\author{Eric B. Kahn}

\begin{document}
\begin{abstract}
 We give a conjectural description for the kernel of the map assigning to each finite $\mathbb Z_p$-free $G\times\mathbb Z_p$-set its rational permutation module where $G$ is a finite $p$-group. We prove that this conjecture is true when $G$ is an elementary abelian $p$-group or a cyclic $p$-group.
\end{abstract}

\maketitle

\section{Introduction}

Given a prime $p$ and a finite $p$-group $G$, there is a natural map between the Burnside ring $A(G)$ and the rational representation ring $R(G)$ taking finite $G$-sets $X$ to their permutation modules $\mathbb Q[X]$. It was shown by J. Ritter and G. Segal \cite{JR, GS} that when $G$ is a $p$-group this map is a surjection and we can study the connection between $A(G)$ and $R(G)$ by looking at generators for the kernel denoted $N(G)$ of the map. These generators were shown by Tornehave in \cite{JT} to always be induced by ``small'' subquotients of $G$.

From a topological viewpoint we have the Segal conjecture proved by G. Carlsson in \cite{GC} which states that the ring $A(G)$ is isomorphic up to completion to the stable cohomotopy group of $BG$. In addition it was shown by M. Atiyah in \cite{MA} that the complex representation ring $R_{\mathbb C}(G)$ is isomorphic up to completion to the complex $K$-theory of the classifying space $BG$. Thus we gain a connection between stable cohomotopy theory and complex $K$-theory of $BG$ which can be studied by purely algebraically means.

In a generalization of these ideas, we can consider for a finite $p$-group $G$, the natural map between the relative Burnside ring $A(G,\mathbb Z_p)$ and the relative rational representation ring $R(G,\mathbb Z_p)$ which takes finite, $\mathbb Z_p$-free, $G\times \mathbb Z_p$-sets $X$ to the $\mathbb Q[\mathbb Z_p]$-free, $G\times \mathbb Z_p$ permutation $\mathbb Q[G\times \mathbb Z_p]$-modules $\mathbb Q[X]$. We give their definition in section~\ref{sec:rel}. It was show by M. Anton in \cite{MFA} that this map is surjective.

The purpose of this paper is to give a conjecture describing the kernel denoted $N(G,\mathbb Z_p)$ of the map in the relative case and to prove the conjecture for some specific $p$-groups $G$. We begin in section~\ref{sec:prelim} by giving a brief overview of the classical case, including the construction of induction maps. In section~\ref{sec:elemabel} we consider the case when $G$ is an elementary abelian $p$-group. In this case we decompose both the classical and relative Burnside rings into graded modules and compute their ranks. The main theorem is given in section ~\ref{sec:mainthrm} and describes the generators of the relative kernel $N(G,\mathbb Z_p)$. In section ~\ref{sec:final} we offer a conjecture for describing the kernel of the relative map for any $p$-group $G$ and offer proofs in the cases of $G$ a cyclic $p$-group and an elementary abelian $p$-group.

\textbf{Acknowledgments} I am grateful to the department at the University of Kentucky for funding over the summer to work on this problem and to Dr. M. Anton for his enlightening discussions.

\section{Preliminaries}\label{sec:prelim}

\subsection{Burnside and representation rings}

For a finite group $G$ the isomorphism classes of finite $G$-sets form a semiring $S$ with respect to disjoint union and direct product. The Burnside ring $A(G)$ is defined to be the Grothendieck construction of the semiring $S$. In fact $A(G)$ is a free $\mathbb Z$-module with a basis given by the set of left coset spaces $[G/L]$ where $L$ runs thru conjugacy class representatives of subgroups $L<G$. For each such subgroup we define an induction map $L\uparrow:A(L)\to A(G)$ by sending an $L$-set $X$ to the $G$-set $G\times_LX$ where $gl\times x=g\times lx$ for all $(g,l,x)$ in $G\times L\times X$. This definition extends to induction maps $L/C\uparrow:A(L/C)\to A(G)$ via the pullback map $A(L/C)\to A(L)$ where $L/C$ is a subquotient of $G$. The induction maps are $\mathbb Z$-linear but do not preserve the product.

Likewise let $T$ be the semiring of isomorphism classes of finitely generated $\mathbb Q[G]$-modules with respect to direct sum and tensor product. The rational representation ring $R(G)$ is defined to be the Grothendieck construction of $T$. For $L<G$, the induction map $L\uparrow:R(L)\to R(G)$ sends a $\mathbb Q[L]$-module $M$ to the $\mathbb Q[G]$-module $\mathbb Q[G]\otimes_{\mathbb Q[L]}M$ and extends to subquotients as in the Burnside ring case.

The Burnside and representation rings are related by a natural ring homomorphism $f:A(G)\to R(G)$ sending a $G$-set $X$ to the permutation $\mathbb Q[G]$-module $\mathbb Q[X]$. It is immediate that $f$ commutes with the induction maps.

\begin{definition}
The \emph{Burnside kernel} $N(G)$ is the kernel of the map $f$.
\end{definition}

\subsection{Relative Burnside and representation modules}\label{sec:rel}

If $\tilde G=G\times H$ is a direct product of two finite groups then a $\tilde G$-set is thought of with $G$ acting on the left and $H$ on the right.  Let $S'$ be the monoid of isomorphism classes of finite $H$-free $\tilde G$-sets with respect to disjoint union. The relative Burnside module $A(G,H)$ is the Grothendieck construction of the monoid $S'$. Then $A(G,H)\subset A(\tilde G)$ is a free $\mathbb Z$-submodule with a basis given by twisted products $[G\times_\rho H]$ where $\rho$ runs thru conjugacy class representatives of homomorphisms $\rho:K\to H$ with $K<G$ and $gk\times h=g\times\rho(k)h$ for all $(g,k,h)$ in $G\times K\times H$.

Similarly, a $\mathbb Q[\tilde G]$-module is thought of with $\mathbb Q[G]$ acting on the left and $\mathbb Q[H]$ on the right. Let $T'$ be the monoid of isomorphism classes of finitely generated $\mathbb Q[H]$-free $\mathbb Q[\tilde G]$-modules with respect to direct sum. The relative rational representation module $R(G,H)$ is the Grothendieck construction of the monoid $T'$. Then the natural ring homomorphism $f:A(\tilde G)\to R(\tilde G)$ will restrict to a module homomorphism $f':A(G,H)\to R(G,H)$.

\begin{definition} The \emph{relative Burnside kernel} $N(G,H)$ is the kernel of $f'$.\end{definition}

\subsection{Relative induction}\label{sec:relind}

The relative induction maps $L/C\uparrow:\tilde A(L/C)\to A(G,H)$ are defined by the usual induction $L/C\uparrow$ restricted to the submodule $\tilde A(L/C)$ made of those elements of $A(L/C)$ that land in $A(G,H)$ where $L/C$ is a subquotient of $\tilde G$. The same observation applies to the relative induction map $L/C\uparrow:\tilde R(L/C)\to R(G,H)$. It is immediate that the natural map $f'$ from the relative Burnside module to the relative representation module commutes with the relative induction maps. We conclude by proving from scratch the following

\begin{proposition}
All induction maps are injective.
\end{proposition}

\begin{proof}
 Let $M$ be any of the monoids $S$, $T$, $S'$, $T'$ defining the Burnside and representation modules and their relative versions. The induction map is defined by a homomorphism $L/C\uparrow:M\to N$ of monoids extended to the Grothendieck constructions where $L/C$ is a subquotient of $\tilde G$. The Grothendieck construction of $M$ consists of fomal differences $[X]-[Y]$ of elements in $M$ such that $[X]-[Y]=[X']-[Y']$ if and only if $$[X+Y'+Z]=[X'+Y+Z]$$ for some $[Z]$ in $M$. In particular, if $L/C\uparrow [X]-L/C\uparrow[Y]=0$ then
$$[L/C\uparrow X+V]=[L/C\uparrow Y+V]$$
for some $[V]$ in $N$. By restricting the $G$-structure to an $L$-structure we have a restriction map $L\downarrow$ such that $L\downarrow (L/C\uparrow [X])=[\tilde G:L][X]$. In particular,
$$
[\tilde G:L][X]+L\downarrow[V]=[\tilde G:L][Y]+L\downarrow[V].
$$
Since each element of $M$ has a unique decomposition into a sum of irreducible elements, we conclude that $[X]=[Y]$ proving the injectivity of the induction map.
\end{proof}

\subsection{The $p$-group case}\label{i}

For $G$ a finite $p$-group  it was shown by Tornehave \cite{JT} that $N(G)$ is generated by the induced kernels $L/C\uparrow N(L/C)$ where $L/C$ runs thru all subquotients of $G$ that are isomorphic to the elementary abelian group $\mathbb Z_p\times\mathbb Z_p$, the dihedral group, or the nonabelian group of order $p^3$ and exponent $p$. Combining this with the Ritter-Segal \cite{JR,GS} proof for the surjectivity of $f$ we get a well understood short exact sequence:
\begin{equation}\label{s}
0\to N(G)\to A(G)\xrightarrow{f} R(G)\to 0.
\end{equation}
In the abelian case $G=\mathbb Z_p\times\mathbb Z_p$ for instance, it is shown in \cite{EL} that $N(G)$ is the free cyclic group generated by
\begin{equation}\label{t}
[G]-\sum[G/C]+p[G/G]
\end{equation}
where $C$ runs thru all proper cyclic subgroups of $G$. In the relative case $\tilde G=G\times H$ with $G$ a finite $p$-group and $H=\mathbb Z_p$ it is known by \cite{MFA} only that we have a short exact sequence
\begin{equation}\label{a}
0\to N(G,H)\to A(G,H)\xrightarrow{f'} R(G,H)\to 0
\end{equation}
and the purpose of this paper is to study $N(G,H)$.

\subsection{A useful trick}

\begin{lemma}\label{lem:exact}
Consider the chain complex of finitely generated free $\mathbb Z$-modules
\[0\stackrel{}{\rightarrow}A\stackrel{\alpha}{\rightarrow}B\stackrel{\beta}{\rightarrow}C\stackrel{}{\rightarrow}0 \]
with $\alpha$ injective and $\beta$ surjective. If the cokernel of $\alpha$ is a free module and the rank of the image of $\alpha$ equals the rank of the kernel of $\beta$, then the sequence is exact.
\end{lemma}

\begin{proof}
Since $Im(\alpha)\subset Ker(\beta)$ and $Coker(\alpha)$ is free, we have the free $\mathbb Z$-submodule $Ker(\beta)/Im(\alpha)\subset B/Im(\alpha)$. But the rank of the image of $\alpha$ equals the rank of the kernel of $\beta$ so that $Ker(\beta)/Im(\alpha)$ is torsion. Therefore $Ker(\beta)/Im(\alpha)=0$.
\end{proof}

\section{The Relative Burnside Kernel for Elementary Abelian Groups}\label{sec:elemabel}

\subsection{Notations}\label{sec:grade}
In this section let $G=\mathbb Z_p^n$ and $H=\mathbb Z_p$ so that $\tilde G=G\times H=\mathbb Z_p^{n+1}$. Also, we denote
$$
A=A(\tilde G),\ R=R(\tilde G),\ N=N(\tilde G)
$$
$$
A'=A(G,H),\ R'=R(G,H),\ N'=N(G,H)
$$
and define $A_k\subset A$ to be generated by all $[\tilde G/L]$ with $L\subset\tilde G$ of dimension $k$. Thus,
$$
A=A_0\oplus A_1\oplus...\oplus A_{n+1}
$$
and a similar decomposition holds for $A'$ with $A'_k=A_k\cap A'$.

\subsection{Rank calculations}\label{sec:rank}

Let $G(k,n)$ denote the number of $k$-dimensional subspaces of the vector space $\mathbb Z_p^n$. Then by \cite[p. 28]{RPS} we have
$$
G(k,n)=\prod_{j=1}^k\frac{p^{n-j+1}-1}{p^j-1}
$$

\begin{proposition}\label{ak} The ranks $a_k$ and $a'_k$ of $A_k$ and $A'_k$ are given by the formulas
$$
a_k=G(k,n+1),\ a'_k=p^kG(k,n).
$$
\end{proposition}

\begin{proof}
The basis elements $[\tilde G/L]$ for $A_k$ are in one-to-one correspondence with the $k$-dimensional subspaces $L<\tilde G=\mathbb Z_p^{n+1}$. Hence, we get the first formula.

The basis elements $[G\times_\rho H]$ for $A'_k$ are in one-to-one correspondence with pairs $(K,\rho)$ with $K< G$ a $k$-dimensional subspace and $\rho:K\to H$ a homomorphism. Given $K$, $\rho$ is uniquely determined by its kernel and an automorphism of its image. If $K$ is $k$-dimensional, the kernel of $\rho$ is either $K$ or any $(k-1)$-dimensional subspace of $K$. In the later case the image admits $(p-1)$ automorphisms. Hence, for a given $k$-dimensional $K$ there are $(p-1)G(k-1,k)+1$ different $\rho$'s. For a given dimension $k$ the number of pairs $(K,\rho)$ is thus given by the formula:
$$
G(k,n)[(p-1)G(k-1,k)+1]=G(k,n)[(p-1)\sum_{i=0}^{k-1}p^i+1]=p^kG(k,n).
$$
\end{proof}

Let $\zeta$ denote a primitive $p$-root of unity and $F=\mathbb Q[\zeta]$ be the associated cyclotomic field. For each $s\in\mathbb Z_p^{n+1}$  let $F_s$ be the $\mathbb Q[\tilde G]$-module $F$ obtained by letting the $i^{th}$ canonical generator of $\tilde G$ act on $F$ via the automorphism sending $\zeta$ to $\zeta^{s_i}$ where $s_i$ is the $i^{th}$ coordinate of $s$.

\begin{proposition}\label{r}
The ranks $r$ and $r'$ of $R$ and $R'$ are given by the formulas
$$
r=G(1,n+1)+1,\ r'=G(1,n+1).
$$
\end{proposition}
\begin{proof}
With the above notations, two isomorphism classes are equal $[F_s]=[F_t]$ if and only if $t=s=0$ or $t=us$ for some unit $u$ in $\mathbb Z_p$. In the later case we say that $s$ and $t$ represent the same point $[s]=[t]$ in the projective $n$-space $P^n$ over $\mathbb Z_p$. With this observation $[F_s]$ indexed by $[s]\in P^n$ and the trivial module $[\mathbb Q]$ form a basis for $R$. Thus we get the first formula.

For the second formula we claim that a basis for $R'$ is given by the elements
$$
[F_{s'\times 1}]+[\mathbb Q],\ \ [F_{t\times 0}]-(p-1)[\mathbb Q]
$$
indexed by $s'\in\mathbb Z_p^n$ and $[t]\in P^{n-1}$. Let $\mathcal{B}$ denote the set of these elements and $\mathcal M$ the $\mathbb Z$-module generated by $\mathcal B$. Since $F_{0\times 1}+\mathbb Q=\mathbb Q[H]$ it follows that by forgetting the $G$-action, the elements:
\[ [F_{s'\times 1}]+[\mathbb Q],\,\, [F_{t\times 0}]+(p-1)[F_{0\times 1}],\text{ and } (p-1)[F_{0\times 1}]+(p-1)[\mathbb Q] \]
are all represented by the $\mathbb Q[H]$-free modules $\mathbb Q[H]$ or $(p-1)\mathbb Q[H]$. Thus $\mathcal M\subset R'$ and it is immediate that $\mathcal B$ is a linearly independent set. Now by inspection $R/\mathcal M$ is the free module generated by $[\mathbb Q]$ and $m[\mathbb Q]\in R'$ implies $m=0$. Thus the rank of $R'$ equals the rank of $\mathcal M$. In particular Lemma~\ref{lem:exact} applies to the sequence:
\[0\rightarrow\mathcal M {\rightarrow}R{\rightarrow}R/R'\rightarrow0 \]
implying that $R'=\mathcal M$.
\end{proof}

From the Propositions \ref{ak} and \ref{r} and the short exact sequences \eqref{s} and \eqref{a} we deduce the following result.
\begin{corollary}\label{n}
The ranks $b$ and $b'$ of $N$ and $N'$ are given by the formulas
$$
b=\sum_{k=0}^{n-1}G(k,n+1),\ b'=\sum_{k=0}^np^kG(k,n)-G(1,n+1).
$$
\end{corollary}

\subsection{The main theorem}\label{sec:mainthrm}

It is convenient to identify each basis element $[\tilde G/L]$ of $A$ where $L<\tilde G$ with the projective subspace $(L)\subset P^n$ generated by $L$. Also, let $e$ denote the distinguished vector $(0,...,0,1)\in\mathbb Z_p^{n+1}$. Then we have the following characterization for the basis elements of $A'_k$:

\begin{lemma}\label{lem:P^n}
 The submodule $A'_k\subset A'$ is the free abelian group on the set of projective subspaces $(L)\subset P^n$ with $L<\tilde G$ of dimension $k$ not containing $e$.
\end{lemma}

\begin{proof} It is easy to see that the basis elements $[G\times_\rho H]$ of $A'$ associated with a pair $(K,\rho)$ is of the form $[\tilde G/L]$ where $K<G$, $\rho:K\to H$ is a homomorphism, and $L=\{(k,\rho(k))|k\in K\}$ is a linear subspace of $\tilde G$ not containing $e$.

Conversely, let $(L)\subset P^n$ with $L<\tilde G$ of dimension $k$ not containing $e$ and define $K$ to be the image of the canonical projection $\tilde G\rightarrow G$. If $(g,h)$ is an element in $L$ which maps to 0 under the projection, then $g=0$. This would imply $he\in L$ so $h=0$. Thus the projection induces an isomorphism $L\cong K$. Let $\alpha:K\rightarrow L$ be the inverse and define $\rho :K\rightarrow H$ by composing $\alpha$ with the canonical projection $\tilde G\rightarrow H$. We can then check that $[\tilde G/L]=[G\times_{\rho} H]$.
\end{proof}

Given $L$ a subspace of codimension at least 2 in $\tilde G$ we define $L^*<\tilde G$
to be a distinguished subspace such that the following two conditions are both
satisfied:
\begin{enumerate}
\item\label{cond1} $L^*$ contains $L$ and $L^*/L$ has rank 2
\item\label{cond2} If $L$ does not contain $e$ and has codimension at least 3 then $L^*$ does not
contain $e$
\end{enumerate}
Now we observe that $L^*$ always exists subject to the two conditions. In
particular, if $L$ has codimension exactly 2 then $L^*=\tilde G$ is the only choice without violating condition 2.

\begin{definition}\label{def:LC}
For each such $L$ define
$$
t(L)=(L)-\sum(C)+p(L^*)
$$
where the sum is over all proper subspaces $L<C<L^*$.
\end{definition}

In particular, define $M_{n-1}$ to be the set of all $(L)$ with $L<\tilde G$ an $(n-1)$-dimensional subpace where $e\not\in L$. By Lemma~\ref{lem:P^n}, this set is also a basis for $A'_{n-1}.$

\begin{definition}\label{def:rel}
Let $A_{n-1}''$ be the submodule of $A'_{n-1}$ generated by all those differences
$(L)-(L')$ of elements in $M_{n-1}$ that are  subject to the relation
\[ (L+\mathbb Z_p e) = (L' + \mathbb Z_p e). \]
\end{definition}

\begin{theorem}\label{ther:ses}
The rank of $A_{n-1}''$ is $G(1,n)(p^{n-1}-1)$ and we have the following commutative diagram of short exact sequences:
$$
\begin{CD}
0@>>>A_0'\oplus A'_1\oplus...\oplus A'_{n-2}\oplus A_{n-1}''@>t'>> A'@>f'>> R' @>>>0\\
@. @VVV @VVV @VVV\\
0@>>> A_0\oplus A_1\oplus...\oplus A_{n-2}\oplus A_{n-1}@>t>> A @>f>> R @>>> 0
\end{CD}
$$
where the vertical arrows are all inclusions and $t'$, $f'$ are the restrictions of $t$, $f$.
\end{theorem}

\begin{proof}
Given $(L)$ in $M_{n-1}$,
\[ t(L)=(L)-(C_0)-\sum (C) +p(\tilde G)\]
where $C_0=L+\mathbb Z_p e$ and the sum is taken over all $L<C<\tilde G$ not containing $e$. By Definition~\ref{def:rel}, if $(L)-(L')$ is a generator of $A_{n-1}''$ then
\[ (L+\mathbb Z_p e) = (L'+\mathbb Z_p e)=(C_0)\]
and we deduce that $t((L)-(L'))$ is in $A'$ so all maps in the diagram are well defined. Also it was shown respectively in \cite{JR, GS} and \cite{MFA} that $f$ and $f'$ are surjective.

From subsection~\ref{i} we know that the kernel of $f$ is generated by the induced kernels $L/C\uparrow N(L/C)$ where $L/C\cong \mathbb Z_p\times \mathbb Z_p$. In particular, by applying $L/C\uparrow$ to equation~\eqref{t} of Section~\ref{i} with $G=L/C$ we deduce that $N(\tilde G)$ is generated by elements of the form
\[ [\tilde G/C]-\sum [\tilde G/D] +p[\tilde G/L]\]
where $L/C$ is any subquotient of $\tilde G$ isomorphic to $\mathbb Z_p\times \mathbb Z_p$ and the sum runs over all proper subgroups $C<D<L$. By Definition~\ref{def:LC} the above elements with $L/C$ replaced by $L^*/L$ generate the image of $t$ so that the composition $f\circ t$ and $f'\circ t'$ are both zero.

We are left to prove the injectivity of the map $t$ and the inclusions of the kernels of $f$ and $f'$ inside the images of $t$ and $t'$ respectively. Under the map $t$, each basis element $(L)$ of $A_i$ is mapped to an element inside $A_i\oplus A_{i+1}\oplus A_{i+2}$ whose first component is again $(L)$. Therefore the matrix representation of $t$ is upper triangular with cokernel $A_n\oplus A_{n+1}$ which is free. Hence $t$, and therefore $t'$, are injective.

Regarding the exactness at $A$ observe that by Proposition~\ref{ak} and the injectivity of $t$ it follows that the rank of $t$ is the sum $G(k,n+1)$ for $k=0,1,...,n-1$. The same sum by Corollary~\ref{n} is the rank of the kernel of $f$. Since the cokernel of $t$ is a free module we conclude by Lemma~\ref{lem:exact} that the bottom sequence is exact at $A$.

To determine exactness at $A'$ we must first determine the rank of $A_{n-1}'$. As calculated in Proposition~\ref{ak}, the rank of $A_{n-1}'$ which is also the order of $M_{n-1}$ is equal to $a_{n-1}'=p^{n-1}G(n-1,n)$. We observe that $M_{n-1}$ breaks into $G(n-1,n)$ equivalence classes relative to the equivalence relation $(L)\sim (L')$ if and only if \[ L+\mathbb Z_pe=L'+\mathbb Z_pe.\]
By Lemma~\ref{lem:P^n} and Proposition~\ref{ak} with $n$ replaced by $n-1$, each $n$-subspace containing $e$ contains $p^{n-1}G(n-1,n-1)$ subspaces of dimension $n-1$ not containing $e$. This is the number of elements in any of the equivalence classes. Hence, since each equivalence class produces $p^{n-1}-1$ basis elements for $A_{n-1}''$ and there are $a_{n-1}'p^{n-1}$ equivalence classes, we conclude that the rank of $A_{n-1}''$ is given by the formula
\[ a_{n-1}'p^{n-1}(p^{n-1}-1)=G(1,n)(p^{n-1}-1).\]
Combining this with the fact that $t'$ is injective it follows that the image of $t'$ has rank equal to the kernel of $f'$. Moreover when considering the generators of $A_{n-1}''$, if we allow any given basis element $(K)\in A(n-1)$ to play the role of an $(L)$ in the difference $(L)-(L')$ at most once, then we see that the matrix of $t'$ will be upper triangular as $t'$ maps a difference $(L)-(L')$ to an element inside $A_{n-1}\oplus A_n\oplus A_{n+1}$ with first component $(L)-(L')$. Therefore the cokernel of $t'$ is a free module and by Lemma~\ref{lem:exact}, the top row is exact.

%
%
%
%
\end{proof}

\subsection{An illustration for $n=2$ and $p=2$}\label{sec:n=p=2}
Order $\mathbb Z_2$ such that $0<1$ and order $\mathbb Z_2^3$ lexicographically. Then for $n=p=2$ we gain a labeling of the basis of $A(\tilde G$) $\{e_i\}$ such that $e_1<e_2<...<e_{16}$. With this labeling of the basis of $A(\tilde G)$, the subgroup lattice of $\tilde G$ can be represented by the graph $E$ below and offers a visual description of the relationship between basis elements $e_i$ and $e_j$.

\begin{center}
\includegraphics[scale=.6]{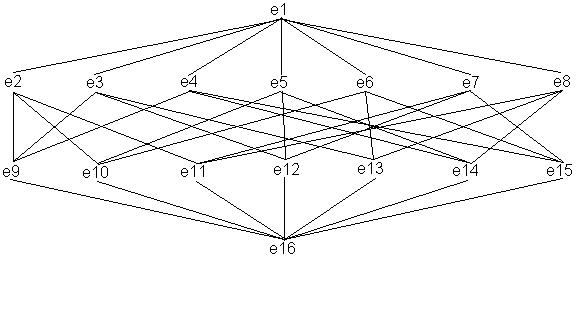}
\end{center}

Theorem~\ref{ther:ses} implies that we have the commutative diagram of short exact sequences:
$$
\begin{CD}
0@>>>A_0'\oplus A_{1}''@>t'>> A'@>f'>> R' @>>>0\\
@. @VVV @VVV @VVV\\
0@>>> A_0\oplus A_1@>t>> A @>f>> R @>>> 0
\end{CD}
$$
We see using our basis that
\[A_0=A_0'=\mathbb Ze_1,\, A_1=\sum_{i=2}^8 \mathbb Ze_i \text{ and } A_1'=\mathbb Z(e_3-e_4) + \mathbb Z(e_5-e_6) + \mathbb Z(e_7-e_8).\]
Hence $t$ is well defined on $A_1$ while we define:
\[ t(e_1)=e_1-e_3-e_5-e_7+2e_{12}.\]
Define the subgraph $E_i$ to be the full subgraph of $E$ where the vertices are the terms occurring in $t(e_i)$. Then the image $t'(e_i-e_j)$ is associated with the subgraph $E_i-E_j$ whose vertices are those in $E_i$ and $E_j$. For example, if $i=3$ and $j=4$ these subgraphs are:

\begin{center}
\includegraphics[scale=.6]{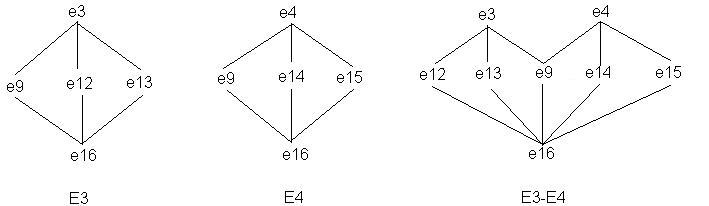}
\end{center}

Conversely, given a subgraph $E_i$ the image $t(e_i)$ is  uniquely determined by taking a weighted sum of the vertices of $E_i$. Moreover, given a subgraph $E_i-E_j$, the image $t(e_i-e_j)$ is also uniquely determined by the vertices of $E_i-E_j$. It follows that the kernel of $f$ is generated by all of the subgraphs $E_i$ for $i=1,2,...,8$ and the kernel of $f'$ is generated by all the non-singular subgraphs $E_1,\, E_3-E_4,\, E_5-E_6,\, E_7-E_8$.

\section{Final Remarks}\label{sec:final}
We would like to develop a description for the kernel $N(G,H)$ with $H\cong \mathbb Z_p$ similar to that given by Tornehave in \cite{JT} for $N(G)$ with $G$ an arbitrary finite $p$-groups $G$. Define $\tilde N(L/C)$ for $L/C$ a subquotient of $\tilde G=G\times H$ to be the intersection of $N(L/C)$ with the submodule $\tilde A(L/C)$ of $A(L/C)$ that lands inside $A(G,H)$ under the induction $L/C\uparrow$ of Section~\ref{sec:relind}.

\begin{conjecture}\label{conj}
Let $p$ be a prime, $G$ any finite $p$-group, and $H\cong \mathbb Z_p$. Then
\[ N(G,H) = \sum L/C\uparrow \tilde N(L/C)\]
where the sum is taken over subquotients $L/C$  of $\tilde G$ isomorphic to $T\times H$ where $T$ is the elementary abelian group $\mathbb Z_p\times \mathbb Z_p$, the dihedral group, or the nonabelian group of order $p^3$ and exponent $p$.
\end{conjecture}

For $G$ elementary abelian or cyclic this conjecture can readily be checked using Theorem~\ref{ther:ses} and rank arguments.

\begin{proposition}
Let $p$ be any prime, $G$ be an elementary abelian $p$-group, and $H\cong \mathbb{Z}_p$. Then
\[ N(G,H) = \sum L/C\uparrow \tilde N(L/C)\]
with the sum taken over all subquotients $L/C\cong \mathbb Z_p^3$.
\end{proposition}

\begin{proof}
From Theorem~\ref{ther:ses} we know that the image of $t$ generates the kernel $N(G,H)$. If $(L)\in A_i'$ with $0\leq i\leq n-2$, then there exists subgroups $L<L^*<B<\tilde G$ such that $B/L\cong \mathbb Z_p^3$ where $L^*$ is the distinguished element used to define $t$ in Definition~\ref{def:LC}. In addition, regardless of our choice of $B$,
\[t((L))\in B/L\uparrow \tilde N(B/L).\]
If $(L)-(L')\in \tilde A_{n-1}$, let
\[C=L+\mathbb Z_pe=L'+\mathbb Z_pe,\, D=L\cap L'.\]
We see immediately that $\tilde G/D\cong \mathbb Z_p^3$ and also that  $t((L)-(L'))$ is an element of $\tilde G/D\uparrow \tilde N(\tilde G/D)$. Hence we conclude that
\[ N(G,H) \subset \sum L/C\uparrow \tilde N(L/C).\]

The converse is immediate.
\end{proof}

\begin{proposition}\label{lem:orderp^k}
Let $p$ be any prime, $G$ the cyclic $p$-group with order $p^k$, and $H\cong \mathbb{Z}_p$. Then $f'$ is an isomorphism between $A(G,H)$ and $R(G,H)$.
\end{proposition}

\begin{proof}
Let $\tilde G = G\times H$ and $\xi$ be the primitive $p^k$-root of unity. Since $G$ is cyclic, easily the rank of $A(G,H)$ is equal to $kp+1$ as $G$ has $k+1$ subgroups and for a nontrivial subgroup $K<G$, there are $p$ homomorphisms $\rho :K\rightarrow H$. Let $F_{\nu,\phi}=\mathbb Q(\xi^{p^{k-\nu}}, \xi^{\phi p^{k-1}})$ be the $\mathbb Q[\tilde G]$-module with the generators of $G$ and $H$ acting by multiplication by $\xi^{p^{k-\nu}}$ and $\xi^{\phi p^{k-1}}$ respectively where $\nu = 0,1,...,k$ and $\phi =0,1,...,p-1$. Then the irreducible $\mathbb Q[\tilde G]$-modules as seen from the decomposition of the group ring $\mathbb Q[\tilde G]$ are:
\begin{eqnarray*}
F_{0,0} & = & \mathbb Q\\
F_{0,1} & = & \mathbb Q(\xi^{\phi p^{k-1}})\\
F_{\nu,\phi} & \text{with} & \text{ $\nu =1,...,k$ and $\phi = 0,1,...,p-1$}.
\end{eqnarray*}
This implies the rank of $R(\tilde G)$ is $kp+2$.

For $[M],[M']\in R(\tilde G)$, define $[M]\equiv [M']$ if we have  $[M]-[M']\in R(G,H)$. Using this relation we immediately gain the following equivalences from \cite{MFA}:
\label{lem:cycsim}
\begin{eqnarray*}
[F_{\nu,\phi}]  & \equiv & -p^{\nu -1} [\mathbb Q] \text{ \hspace{2em} for $\nu=1,...,k$ and $\phi=1,...,p-1$},\\
{[F_{0,1}]}     & \equiv & -[\mathbb Q],\\
{[F_{\nu,0}]}   & \equiv & p^{\nu -1}(p-1) [\mathbb Q] \text{ \hspace{2em} for $\nu=1,...,k$.}
\end{eqnarray*}

The equivalences~\label{lem:cycsim} imply that the rank of $R(\tilde G)/R(G,H)$ is less than or equal to 1. In addition, since $f'$ is surjective, the rank of $A(G,H)$ is $kp+1$, and the rank of $R(\tilde G)$ is $kp+2$, we see that the rank of $R(\tilde G)/R(G,H)$ is at least 1. Thus the rank of $R(\tilde G)/R(G,H)$ is exactly 1 which implies the rank of $R(G,H)=kp+1$. As $A(G,H)$ is a free module, the rank of $A(G,H)$ is equal to the rank of $R(G,H)$, and $f'$ is a surjection, we see that $f'$ is an isomorphism.

\end{proof}
%
%
%
%
%

As a corollary, Conjecture~\ref{conj} is true for $G$ a cyclic $p$-group.

\bibliographystyle{plain}
\bibliography{bib}

\end{document}